\documentclass[12pt]{article}
\usepackage{amsmath}
\usepackage[usenames]{color}
\usepackage{mathrsfs}
\usepackage{amsfonts}
\usepackage{amssymb,amsmath}
\usepackage{CJK}
\usepackage{cite}
\usepackage{cases}
\usepackage{amsthm}
\usepackage{graphicx}
\usepackage{verbatim}
\usepackage{textcomp}
\usepackage{amssymb}
\usepackage{cite}
\usepackage[section]{placeins}
\usepackage{amsmath}
\usepackage{latexsym}
\usepackage{amscd}
\usepackage{amsthm}
\usepackage{mathrsfs}
\usepackage{xypic}
\usepackage{bm}
\usepackage{url}
\usepackage{hyperref}
\usepackage{multirow}
\def\cl{\centerline}

\def\vs{\vspace*}

\def\ni{\noindent}

\numberwithin{equation}{section}
\newtheorem{theo}{Theorem}[section]

\newtheorem{exam}[theo]{Example}

\newtheorem{remark}[theo]{Remark}
\newtheorem{thm}{Theorem}[section]
\newtheorem{corr}{Corollary}[section]
\newtheorem{lem}{Lemma}[section]

\theoremstyle{definition}
\newtheorem{defn}{Definition}[section]

\theoremstyle{remark}

\numberwithin{equation}{section} 

\begin{document}
\begin{center}
\cl{\large\bf \vs{6pt}  Isoparametric hypersurfaces in conic Finsler manifolds}
\footnote {\ \ $^\dag\,$ The corresponding author's email: dpl2021@163.com
}
\cl{ Qun He$^{1}$, Xin Huang$^{1}$, Peilong Dong$^{2,\dag}$}

\cl{\small 1. School of Mathematical Sciences, Tongji University, Shanghai
200092, China}
\cl{\small 2. School of Mathematics and Statistics, Zhengzhou Normal University,} \cl{\small Zhengzhou Henan 450044, China.}
\end{center}

{\small
\parskip .005 truein
\baselineskip 3pt \lineskip 3pt

\noindent{{\bf Abstract:}
In this paper, we introduce isoparametric functions and isoparametric hypersurfaces in conic Finsler spaces.We find that there are probably other isoparametric hypersurfaces in conic Minkowski spaces besides the conic Minkowski hyperplanes, conic Minkowski hyperspheres and conic Minkowski cylinders, such as helicoids. Moreover, we give a complete classification of isoparametric hypersurfaces in kropina spaces with constant flag curvature.\\
\vs{5pt}

\ni{\bf Key words}:
Isoparametric hypersurface, conic Finsler manifolds, Kropina space, principal curvature.}

\ni{\it Mathematics Subject Classification (2010):} 53B40, 53B25.}
\parskip .001 truein\baselineskip 6pt \lineskip 6pt
\section{Introduction}
Finsler metrics are generalized Riemannian metrics which have no quadratic restriction. The standard definition of a Finsler metric $F$ on a manifold $M$ entails that $F$ is defined on the whole tangent bundle $TM$ and that strong convexity is satisfied, i.e its fundamental tensor $g$ is positive definite. However, in many cases, the metric $F$ is defined only in some conic domain $AM \subsetneq TM$, metrics of which kind are called conic Finsler metrics. It's well known that Kropina metrics are conic Finsler metrics. In a conic Finsler manifold $(M,F)$, the hypersurfaces whose normal vector belongs to $AM$ are called\emph{conic hypersurfaces}.

In Riemannian geometry, the classification of isoparametric hypersurfaces in space forms is a classical geometric problem with a history of almost one hundred years.\cite{GT,MG,QM2} In Finsler geometry, the conception of isoparametric hypersurfaces has been introduced in\cite{HYS}, and the classifications of isoparametric hypersurfaces in some special Finsler spaces have been obtained\cite{HYS,HD,HYS2}. M. Xu and his collaborators prove that under homothetic navigation transformation, hypersurfaces are locally isoparametric if and only if they are locally isoparametric with respect to the original metric\cite{XM2}. As far as we know, isoparametric hypersurfaces in conic Finsler spaces have not been studied.

In this paper, we study isoparametric functions and isoparametric hypersurfaces in conic Finsler spaces. Firstly, we are concerned with the existence of new isoparametric hypersurfaces in conic Minkowski spaces and get the following theorems.
\begin{thm}\label{011}
In a conic Minkowski space $(V,F)$, conic Minkowski hyperplanes, conic Minkowski hyperspheres and conic Minkowski cylinders must be isoparametric hypersurfaces with one or two distinct constant principal curvatures.
\end{thm}
\begin{remark}
When the number of principal curvatures $\tilde{g}= 1$, theorem\ref{011} is also true on the contrary. When the number of principal curvatures $\tilde{g}= 2$, theorem\ref{011}               is not true on the contrary. Theorem\ref{012} is a counterexample
\end{remark}
\begin{thm}\label{012} Let~$(\mathbb{R}^3,F)$~be a $3$-dimension conic Minkowski-$(\alpha,\beta)$ space and the dual metric of $F$ be $F^{\ast}=\alpha ^{\ast}\phi(\frac{\beta^{\ast}} {\alpha^{\ast}})$, where
\begin{equation}\label{139}
\phi=\Big(\frac{\sqrt{b^2-(1+a^2)s^2}}{b}-\frac{as}{b}\arctan\frac{\sqrt{b^2-(1+a^2)s^2}}{as}\Big),~{0<|s|<\frac{b}{\sqrt{1+a^2}}},
\end{equation}
 $\beta^{\ast}=(0,0,b)$ and $a,b$ are two positive constants.~Then in~$(\mathbb{R}^3,F)$,~helicoid~$\textbf r=(u\cos v,u\sin v,av)$ $(0<u<1)$ is a local minimal isoparametric hypersurface with constant principal curvatures $\pm 1$.
\end{thm}
Since both Kropina metrics and Randers metrics can be characterized as the solutions of the Zermelo navigation problem on a Riemannian space $(M, h)$ under the influence of a vector field $W$, we can also give the classifications of isoparametric hypersurfaces in a Kropina space with constant flag curvature by studying the relationship of their principal curvatures.
\begin{thm}\label{130}
Let $(M,F)$ be a Kropina space of constant flag curvature with the navigation data $(h, W)$, then (BH-)isoparametric hypersurfaces of $(M, F, d\mu_{BH})$ must be isoparametric hypersurfaces of $(M,h)$ whose unit normal vector $\bar{\textbf{n}}\neq-W$ and vice versa. Moreover, the number of distinct principal curvatures and the multiplicities of each principal curvature are the same. So the (BH-)isoparametric hypersurfaces in $(M, F, d\mu_{BH })$ can be completely classified (see Table 1 for the accurate classifications).
\end{thm}
The contents of this paper are organized as follows. In Section 2, some fundamental concepts and formulas are given for later use. In Section 3, we give the definition and general properties of isoparametric hypersurfaces in conic Finsler spaces. In Section 4, we consider the isoparametric hypersurfaces in conic Minkowski spaces, and give a new example of isoparametric hypersurfaces in a conic Minkowski-$(\alpha,\beta)$ space. In Section 5, we consider the principal curvatures of hypersurfaces with respect to $F$ and $h$, and derive the classifications of isoparametric hypersurfaces in Kropina spaces with constant flag curvature.

\section{Preliminaries}
\subsection{Conic Finsler metrics.}
In this section, we briefly recall the fundamentals of conic Finsler geometry by Miguel Angel Javaloyesn.\cite{MAJ}
\begin{defn}
Let $M$ be an $m$-dimensional differentiable manifold and $\tilde{A}M \subset TM$ be an open subset of the tangent bundle $TM$ such that $\pi(\tilde{A}M) = M$, where $\pi : TM \rightarrow M$ is the natural projection, and $\tilde{A}M$ is conic in $TM$, i.e. for each $x\in M$, $\tilde{A}_xM :=\tilde{A}M\cap T_xM$ is a conic domain in $T_xM.$ Assume a continous function $F : \tilde{A}M\rightarrow[0, +\infty)$ satisfies\\
(1) $F$ is smooth on $\tilde{A}M\setminus{0}$.\\
(2) $F(x,\lambda y)=\lambda F(x,y)~for~any ~\lambda >0,~x\in M~and ~ y\in \tilde{A}M$.\\
(3) $g=g_{ij}(x,y)dx^{i} \otimes dx^{j}$ is positive definite  on $\tilde{A}M$, where $g_{ij}(x,y)=\frac{1}{2}[F^{2}] _{y^{i}y^{j}}.$\\
Then $F$ is called a \emph{conic Finsler metric}, $(M^{n},F)$ is called a \emph{conic Finsler manifold} and $g$ is called the \emph{fundamental tensor}.
\end{defn}

\begin{exam}\label{3.0}
Let $\phi:[-b_{0},b_{0}]\setminus E\rightarrow(0,+\infty)$ be a smooth positive function, where $E$ is a closed subset of $[-b_{0},b_{0}]$. Suppose that $\alpha(x,y)=\sqrt{ a_{ij}(x)y^i y^j}$ is a Riemannian metric on an open subset $U\subset R^n$ and $\beta(x,y)=b_i(x)y^i$ is a 1-form satisfying $b:=||\beta||_\alpha\leqslant b_{0}$. Define
\begin{equation}
F(x,y)=\alpha (x,y)\phi(s),   ~~~~s=\frac{\beta(x,y)}{\alpha(x,y)},
\end{equation}
 $\phi$ satisfies the requirements for $\forall s\in [-b,b]\setminus E,\forall b\in[0,b_{0}] $
 \begin{equation}
 \phi(s)-s\phi'(s)>0
 \end{equation}
\begin{equation}
\phi(s)-s\phi'(s)+(b^2-s^2)\phi''(s)>0.
\end{equation}
Then $F$ is called a\emph{conic $(\alpha,\beta)$ metric} with conic domain $$\tilde{A}_xM=\{y\in T_xM~|~\frac{\beta(x, y)}{\alpha(x, y)}\in[-b_{0},b_{0}]\setminus E\}.$$
\end{exam}
Let~$(M,F)$ be an~$m$-dimensional oriented smooth conic Finsler manifold and~$TM$ be the tangent bundle over~$M$ with local coordinates $(x,y)$, where~$x=(x^1,\cdots ,x^m)$ and~$y=(y^1,\cdots,y^m)$. Now we will use the following convention of index ranges unless other states:
$$1\leq i, j,\cdots\leq m ;~~~~~~~1\leq a, b, \cdots\leq n< m.$$\\
$\pi_{\tilde{A}}$ is the restriction to $\tilde{A}M$ of the natural projection $\pi:TM\rightarrow M$,  which gives rise to the pull-back bundle~$\pi_{\tilde{A}}^{\ast}TM$~and its dual bundle~$\pi_{\tilde{A}}^{\ast}T^{\ast} M$~over~$\tilde{A}M$~.As in the classical case,~on the pull-back bundle~$\pi_{\tilde{A}}^{\ast}TM$~there exists uniquely the \emph{Chern connection}~$\nabla$ with~$\nabla\frac{\partial}{\partial x^i}=\omega_{i}^{j}\frac{\partial}{\partial
x^{j}}=\Gamma^{i}_{jk}dx^k\otimes\frac{\partial}{\partial x^{j}}$~satisfying
$$ \omega^i_j \wedge dx^j=0,$$
$$dg_{ij}-g_{ik}\omega^{k}_{j}-g_{kj}\omega^{k}_{i}=2C_{ijk}(dy^{k}+N^{k}_{l}dx^{l}),$$
$$~~~~N_j^i:=\frac{\partial G^i}{\partial y^j}=\Gamma^{i}_{jk}y^k,~~~~~\delta y^{l}:=\frac{1}{F}(dy^l+y^j\omega^{l}_{j}),$$
where~$C_{ijk}=\frac{1}{2}\frac{\partial g_{ij}}{\partial y^k}$~is called the\emph{Cartan tensor} and $$G^{i}=\frac{1}{4}g^{il}\left\{[F^{2}]_{x^{k}y^{l}}y^{k}-[F^{2}]_{x^{l}}\right\}$$ are the \emph{geodesic coefficients} of $(M,F)$. The\emph{curvature 2-forms of the Chern connection $\nabla$}~are
\begin{equation}\label{3}
d\omega_{i}^{j}-\omega_j^k\wedge  \omega_k^i=\Omega_j^i=\frac{1}{2}R^{~i}_{j~kl}dx^{k}\wedge
dx^{l}+P^{~i}_{j~kl}dx^{k}\wedge\delta y^{l},
\end{equation}
where $R^{i}_{j~kl}=-R^{i}_{j~lk}$.

For a fixed point $(x,y)\in \tilde{A}M$, let $\varPi _y (v)=span\{y, v\}\subset T_xM$ be a two-dimensional plane in $T_xM$ with the flagpole $y$ at $x\in M$. The flag curvature of $\varPi _y (v)$ is defined by
$$ K(\varPi _y )(v) :=\frac{-R_{ijkl } y^{i}v^{j}y^{k}v^{l}}{(g_{ik}g_{jl}-g_{il}g_{jk})y^{i}v^{j}y^{k}v^{l}}.$$
$(M,F)$ is said to have constant flag curvature if $K(\varPi _y )(v)$ = constant everywhere.

Let $X=X^{i}\frac{\partial}{\partial x^{i}}$ be a vector field, the\emph{covariant derivative} of $X$ along $v=v^i\frac{\partial}{\partial x^i}\in T_xM $ with respect to~$w\in\tilde{A}_xM$~is defined by
\begin{align}D^{w}_{v}X(x):=\left\{v^{j}\frac{\partial X^{i}}{\partial x^{j}}(x)+\Gamma^{i}_{jk}(w)v^{j}X^{k}(x)\right
\}\frac {\partial}{\partial x^{i}},\label{Z1}\end{align} where $\Gamma^{i}_{jk}$ denote the \emph{connection coefficients} of the Chern connection.

\subsection{Legendre transformation}
The \emph{Legendre transformation} of a conic Finsler metric $(M, F )$ is the map $\mathcal{L}:\tilde{A}M\rightarrow \tilde{A^\ast }M, $satisfying $\mathcal{L}(\lambda y)=\lambda \mathcal{L}(y)$ for all~$\lambda > 0,~y\in \tilde{A}_xM$,~and
\begin{equation}\label{20}
\mathcal{L}:\tilde{A}M\rightarrow {T^\ast M},~~~~~~
\mathcal{L}(y)=F(y)[F]_{y^i}(y)dx^i,~~~~\forall y\in \tilde{A}M.                                                                                                                                                                   \end{equation}
It is notable that~$\mathcal{L}$~may not be injective, so~$\mathcal{L}^{-1}$~is multivalued. In order to define the dual metric~$F^\ast$~of~$F$, we can select a conic open subset~$AM$~of~$\tilde{A}M$~so that is a injective function on~$AM$. Set~$A^\ast M=\mathcal{L}(AM)$, then~$\mathcal{L}: AM\rightarrow A^\ast M$~is a differential homeomorphism. Such a cone is called the \emph{applicable cone}. This paper mainly considers the applicable cone. The dual of the Finsler metric $F$ is the function $F^\ast :A^\ast M\rightarrow[0,+\infty)$ defined by
$$F^\ast=F\circ\mathcal{L}^{-1}.$$
Then
\begin{equation}\label{21}
 \mathcal{L}^{-1}(\xi)=F^\ast(\xi)[F^\ast]_{\xi_i}(\xi)\frac{\partial}{\partial x^i},~~~~\forall \xi\in A_x^\ast M=\mathcal{L}(A_xM).
\end{equation}
For a smooth function $f:M\rightarrow R$, the \emph{conic gradient vector} of $f$ at $x\in M$ is defined as $\nabla f(x):={\mathcal
L}^{-1}(df(x))\in A_{x}M$, which can be written as
\begin{align}\nabla f(x):=\left\{\begin{array}{l}
g^{ij}(x,\nabla f)\frac{\partial f}{\partial x^{j}}\frac{\partial }{\partial x^{i}},~~~~~~ df(x) \neq 0,df\in A_x^\ast M,\\
0,~~~~~~~~~~~~~~~~~~~~~~~~~~~~~~~~df(x)=0.
\end{array}\right.\label{Z2}\end{align}
Set $M_f :=\{x\in M | df(x) \neq 0, df\in A_x^\ast M\}$ and $\nabla^{2}f(x) = D^{ \nabla f} (\nabla f)(x)$~for $x\in M_f$, the \emph{ Finsler-Laplacian of $f$} with respect to the volume form $d\mu=\sigma(x)dx^1\wedge dx^2\wedge...\wedge dx^m$ is defined by
$$\Delta_{\sigma} f =\textmd{div}_\sigma(\nabla f).$$
Another \emph{nonlinear Finsler-Laplacian in $M_f$} is defined as
\begin{equation}
\hat \Delta f : =\textmd{tr}_{g_{_{\nabla f}}}(\nabla^{2}f).
 \end{equation}
\begin{lem}\cite{ZS}\label{010}
$\Delta_{\sigma} f =\hat \Delta f -S(\nabla f)$, where
$$S(x,y)=\frac{\partial G^i}{\partial y^i}-y^i\frac{\partial}{\partial x^i}(\ln\sigma(x)).$$
\end{lem}
\section{Conic hypersurfaces of conic Finsler manifolds}
\subsection{Conic submanifolds and conic hypersurfaces}
Let $(M,F)$ be an $m$-dimensional conic Finsler manifold and $\Phi: N\to (M,F)$ be an $n$-dimensional immersion. For simplicity, we will denote $d\Phi X$ by $X$ locally.  The \emph{conormal bundle} of $N$ is
\begin{equation}\label{ss1}
\mathcal{V}(N)=\{(x,\xi)~|~x\in N,\xi\in T_x^{*}M,\xi (X)=0,~\forall~X\in T_xN\}.
\end{equation}
If $\forall x\in M, A_x^\ast M\cap \mathcal{V}_x(N)\neq\varnothing,$
set $\mathcal{N}(N)={\mathcal L}^{-1}(\mathcal{V}(N)\cap A^*M)$. Moreover, we denote the {\it unit normal bundle} of $N$ by
$$\mathcal{V}^0(N)=\{\nu\in \mathcal{V}(N)\cap A^*M|F^*(\nu)=1\},$$ and let
$\mathcal{N}^{0}(N)={\mathcal L}^{-1}(\mathcal{V}^{0}(N))=\{\textbf{n}~|~\textbf{n}={\mathcal L}^{-1}(\nu),~\nu \in \mathcal{V}^{0}(N)\}.$ We call $\textbf{n}\in \mathcal{N}^{0}(N)$ the\emph{unit normal vector} of $N$ and $(N,g_\textbf{n})$ \emph{conic submanifolds}.\\

For any $X\in T_xN$ and $\textbf{n}$, a local smooth section of $\mathcal{N}^{0}(N)$, the \emph{shape operator} ${\mathcal A}_{\textbf{n}}:T_xN\rightarrow T_xN$ is defined by
\begin{equation}
{\mathcal A}_{\textbf{n}}(X)=-[D_{X}^{\textbf n}\textbf n]_{g_\textbf n}^T.
\end{equation}
We call the eigenvalues of ${\mathcal A}_{\textbf{n}}$, $k_1,k_2,\cdots,k_{n}$, the \emph{principal curvatures} and $\hat H_{\textbf{n}}=k_1+k_2+\cdots+k_{n}$ the \emph{mean curvatures} with respect to $\textbf{n}$.

If $k_1=k_2=\cdots=k_{n}$, we call $N$ \emph{totally umbilic}. If $\hat H_{\textbf{n}}=0, \forall \textbf{n}\in \mathcal{N}^{0}(N)$, we call $N$\emph{minimal}.

If $m-n=1$, there exists at least one global unit normal vector field on $N$. Let $\textbf{n}$ is a given unit normal vector field and $\hat g:=\Phi^*g_{\textbf{n}}$. Then $(N,\hat g)$ is a Riemannian manifold and called a \emph{conic hypersurface} of $(M,F)$. In this case
\begin{equation}\label{30}{\mathcal A}_{\textbf{n}}(X)=-D_{X}^{\textbf n}\textbf n.
\end{equation}

 Let $d\mu_M=\sigma(x)dx^1\wedge\cdots\wedge dx^{m}$ be an arbitrary volume form on $(M,F)$. The induced volume form on $N$ determined by $d\mu_M$ can be defined by
\begin{align}
d\mu_{\textbf{n}}
={\sigma}({\Phi(u)})\Phi^*(i_{\textbf{n}}(dx^1\wedge\cdots\wedge dx^{m})),~~~~~u\in N,\label{Z21}
\end{align} where $i_{\textbf{n}}$ denotes the inner multiplication with respect to $\textbf{n}$.

As similar to\cite{SZ1}, the first variational formula of the induced metric with respect to the induced volume element is  $$\frac{d{\text{Vol}}(t)}{dt}\Bigr|_{t=0}=-\int_{N}\mathcal{H}_{d\mu_{\textbf{n}}}({X})d\mu_{\textbf{n}},$$ $\mathcal{H}_{d\mu_{\textbf{n}}}$ is called the \emph{$d\mu_{\textbf{n}}$-mean curvature form} of $\Phi$ with
respect to $\textbf{n}$.

Define
\begin{align}
{H}_{\textbf{n}}:=\mathcal{H}_{d\mu_{\textbf{n}}}(\textbf{n})
\label{2.11}\end{align} We call ${H}_{\textbf{n}}$ the
\emph{$d\mu_{\textbf{n}}$-mean curvature} of $N$ in $(M,F)$.
\begin{lem}\cite{ZS}
${H}_{\textbf{n}}=\hat{H}_{\textbf{n}}+S(\textbf{n}).$
\end{lem}

\begin{lem}\label{lem1}
 Let $F$ be a conic $(\alpha,\beta)$ metric,  and set
\begin{equation}\label{6}
F^{\ast}(\xi)=\alpha ^{\ast}\phi(\frac{\beta^{\ast}} {\alpha^{\ast}} ),
\end{equation}
where $\alpha^{\ast}$ is the dual metric of $\alpha$, and $\beta^{\ast}$ is the dual vector of $\beta$. Then for a conic submanifold in $(M,F)$, the normal vector \emph{$\textbf{n}$} and \emph{$\bar{\textbf{n}}$} with respect to $F$ and $\alpha$ satisfy\end{lem}
\begin{equation}\label{7}
\textbf{n}=(\phi-s\phi')\bar{\textbf{n}}+\phi'\beta^{\ast},
\end{equation}
where~$s=\frac{\beta^\ast}{\alpha^\ast}=\beta(\bar{\textbf{n}})$.
\proof
Let $F$ be a conic $(\alpha,\beta)$ metric, then $F^{\ast}$ is also a conic $(\alpha,\beta)$ metric and
$$g^{\ast ij}=\rho a^{ij}+\rho_{0}b^{i}b^{j}+\rho_{1}(b^{i}\alpha ^{\ast}_{\xi_{j}}+b^{j}\alpha ^{\ast}_{\xi_{i}})-s\rho_{1}\alpha ^{\ast}_{\xi_{i}}\alpha ^{\ast}_{\xi_{j}},$$
where
$$\rho=\phi(\phi-s\phi'),  ~~~ \rho_{0}=\phi\phi''+\phi'^{2},~~~\rho_{1}=(\phi-s\phi')\phi'-s\phi\phi''.$$
Set $\nu=\mathcal{L}({\textbf{n}})$, then $\nu_{i}=\lambda\bar{\nu_{i}}$, where $\lambda>0$. Then we have $\textbf{n}^{i}g_{ij}=\lambda a_{ij}\bar{\textbf{n}}^{j}$.

Because $1=F^{\ast}(\nu)=F^{\ast}(\lambda \bar{\nu})=\lambda F^{\ast}(\bar{\nu})=\lambda\alpha^{\ast}(\bar{\nu})\phi(\bar{\nu})
=\lambda\phi(\bar{\nu})$, we can obtain $\lambda=\frac{1}{\phi(\bar{\nu})}.$
According to $g^{ij}(\textbf{n})=g^{\ast ij}(\nu)=g^{\ast ij}( \bar{\nu})$, $s=\beta(\bar{\textbf{n}})$, we have
\begin{align*}
g^{ij}(\textbf{n}) a_{jk}\bar{\textbf{n}}^{k}&=g^{\ast ij}(\bar{\nu})a_{jk}\bar{\textbf{n}}^{k}\\
&=[\rho a^{ij}+\rho_{0}b^{i}b^{j}+\rho_{1}(b^{i}\alpha^{\ast}_{\xi_{j}}+b^{j}\alpha ^{\ast}_{\xi_{i}})-s\rho_{1}\alpha ^{\ast}_{\xi_{i}}\alpha ^{\ast}_{\xi_{j}}]a_{jk}\bar{\textbf{n}}^{k}\\
&=\rho\bar{\textbf{n}}^{i}+\rho_{0}sb^{i}+\rho_{1}(b^{i}+s\bar{\textbf{n}}^{i})-s\rho_{1}\bar{\textbf{n}}^{i}\\
&=\rho\bar{\textbf{n}}^{i}+\phi\phi'b^{i}.
\end{align*}
Then
\begin{align*}
\textbf{n}^{i}&=\lambda g^{ij}(\textbf{n})a_{jk}\bar{\textbf{n}}^{k}\\
&=\lambda[\rho\bar{\textbf{n}}^{i}+\phi\phi'b^{i}] \\
&=(\phi-s\phi')\bar{\textbf{n}}^{i}+\phi'b^{i}.
\end{align*}
\endproof

\subsection{Isoparametric hypersurfaces}
\begin{defn}
Let $f$ be a non-constant $C^1$ function defined on a conic Finsler manifold $(M,F)$ such that $df\in A^*M$ when $df\neq0$ and $f$ is smooth in $M_{f}$. Set $J = f(M_ f)$. The function $f$ is
said to be \emph{$d\mu$-isoparametric} (resp. \emph{isoparametric}) on $(M,F,d\mu)$, where $d\mu=\sigma(x)dx$, if there exist a smooth function $a(t)$ and a continuous function $b(t)$ on $J$ such that
\begin{equation}\label{4} \left\{\begin{aligned}
&F(\nabla f)=a(f),\\
&\Delta f=b(f),
\end{aligned}\right.
\end{equation}
hold for $\Delta {f}=\Delta_{\sigma} f$ (resp. $\Delta {f}=\hat\Delta f$) on $M_{f}$. All the regular level surfaces $N_t= f ^{-1}(t)$ form an \emph{($d\mu$-)isoparametric family}, each of which is
called an \emph{($d\mu$-)isoparametric hypersurface }in $(M,F,d \mu )$.
$f$ is said to be \emph{transnormal} if it only satisfies the first equation of \eqref{4}.
\end{defn}
If for any $x\in N$, there exists a neighborhood $U$ of $x$ and an isoparametric function $f$ defined on $U$, such that  $N\cap U$ is a regular level hypersurface of $f$, then $N$ is called a \emph {local isoparametric hypersurface}.\cite{HD}

From Lemma \ref{010}, we know that if $(M,F,d\mu)$ has constant $\mathbf{S}$-curvature, then~$f$~is an isoparametric function if and only if it is $d\mu$-isoparametric. Similar to the  classical case, from \cite{HYS,HD}, we can obtain the following theorems.
\begin{thm}\label{133} On a conic Finsler manifold $(M,F,d\mu)$, a transnormal function $f$ is ($d\mu$-)isoparametric if and only if each regular level hypersurface $N_t$ of $f$ has constant ($ d\mu_{\bf{n}}$-)mean curvature, where ${\bf{n}}=\frac{\nabla  f }{F(\nabla f )}$.
\end{thm}
\begin{thm} Let $(M,F,d\mu)$ be an $m$-dimensional conic Finsler manifold with constant flag curvature, then a transnormal function $f$ is isoparametric if and only if each regular level surface of $f$ has constant principal curvatures.
\end{thm}
\begin{thm}\label{132}
Let $N$ be a connected and oriented hypersurface embedded in a connected conic Finsler manifold with constant flag curvature, then $N$ is locally isoparametric if and only if its principal curvatures are all constant.
\end{thm}
The detailed proofs of Theorem\ref{133} to Theorem\ref{132} can be seen in\cite{HYS,HD}.
\section{Isoparametric Hypersurfaces in conic Minkowski spaces}
\subsection{Proof of Theorem 1.1}
We suppose that $(V, F, d\mu )$ is an $m$-dimensional conic Minkowski space and $d\mu$ is a volume form such that {\bf S}-curvature vanishes. Let $F^\ast$ be the dual metric of $F$, which is also a conic Minkowski metric, and $g^{\ast ij }(\xi  ) = \frac{1 }{2 }[F^{\ast2} ( \xi )]_ {\xi ^{i }\xi ^{j }}$, then \eqref{4} can be written as
 \begin{equation}\label{5} \left\{\begin{aligned}
&F^\ast(df)=a(f),\\
&g^{\ast ij }(df)f_{ij}=b(f),
\end{aligned}\right.
\end{equation}
where $f_{ij}=\frac {\partial^2f}{\partial x^i\partial x^j}, df\in A^\ast V$.\\

Let $\tilde{g}$ be the number of distinct constant principal curvatures on an isoparametric hypersurface. Similar to\cite{HYS}, we can obtain \\

\begin{thm}\label{400}Let $(V,F)$ be an $n$-dimensional conic Minkowski space, and $N$ be a conic hypersurface in $(V,F)$.~Then $N$ is totally umbilic if and only if it is an
isoparametric hypersurface with $\tilde{g}=1$, which holds if and only if $N$ is either a conic hyperplane, a conic Minkowski hypersphere or a reverse conic Minkowski
hypersphere.
\end{thm}

Let $\bar V^*$ be an $n$-dimensional subspace of $V^*$ such that $\bar V^*\cap A^*V\neq\varnothing$  and $\tilde{F}$ be the dual metric of $F^*|_{\bar V^*\cap A^*V}$ in $\bar V$. Then $(\bar{V}, \tilde{F})$ is also a conic Minkowski space. The cylinde $$\Sigma _r=\{x\in V|\tilde {F}(\bar x)=r,~d\tilde {F}(\bar x)\in A^\ast V\}$$
is called the \emph{conic Minkowski cylinder} of radius $r$ in conic Minkowski space $(V,F)$. 

\begin{thm}\label{300}In an $m$-dimensional conic Minkowski space $(V, F)$, conic Minkowski cylinders must be isoparametric hypersurfaces with with $g=2$. The converse is not necessarily true.
\end{thm}
 Theorem\ref{011} follows from Theorem\ref{400} and Theorem\ref{300}. 
\subsection{Proof of Theorem 1.2}
\proof
Let $F$ be a conic Minkowski-$(\alpha,\beta)$ metric whose dual metric is $F^{\ast}=\alpha ^{\ast}\phi(\frac{\beta^{\ast}} {\alpha^{\ast}})$  and
$\textbf r=(u\cos v,u\sin v,av)$ be a helicoid. Then we have
$$\textbf r_{u}=(\cos v,\sin v,0),~\textbf r_{v}=(-u\sin v,u\cos v,a).$$
A unit normal vector field of $\textbf{r}$ with respect to $\alpha$ is given by
$$\bar{\textbf{n}}=\frac{1}{\sqrt{u^{2}+a^{2}}}(a\sin v,-a\cos v,u),$$
so

\begin{equation}\label{8}
\bar{\textbf{n}}_{u}=(u^{2}+a^{2})^{-\frac{3}{2}}(-au\sin v, au\cos v,a^{2})=\mu_{1}\textbf r_{v},
\end{equation}
\begin{equation}
\bar{\textbf{n}}_{v}=\frac{a}{\sqrt{u^{2}+a^{2}}}(\cos v,\sin v,0)=\mu_{2}\textbf r_{u}.
\end{equation}
where $\mu_{1}=\frac{a}{(u^{2}+a^{2})^{\frac{3}{2}}}, ~~\mu_{2}=\frac{a}{\sqrt{u^{2}+a^{2}}}=\mu_{1}G, ~~G=u^{2}+a^{2}$.

Set $W=\beta^*=(b^1,b^2,b^3)$, by \eqref{7}, the unit normal vector field of $\textbf{r}$ with respect to $F$ is given by
$$\textbf{n}=(\phi(s)-s\phi'(s))\bar{\textbf{n}}+\phi'(s)W,$$
where~$s=\beta(\bar{\textbf{n}})$.
We get
µœ\begin{align}\label{9}
\textbf{n}_{a}&=(\phi'b_{i}\bar{\textbf{n}}^{i}_{a}-b_{i}\bar{\textbf{n}}^{i}_{a}\phi'
-s\phi''b_{i}\bar{\textbf{n}}^{i}_{a})\bar{\textbf{n}}
+(\phi-s\phi')\bar{\textbf{n}}_{a}+\phi''b_{i}\bar{\textbf{n}}^{i}_{a}W\\ \nonumber
&=(\phi-s\phi')\bar{\textbf{n}}_{a}+\phi''b_{i}\bar{\textbf{n}}^{i}_{a}(W-s\bar{\textbf{n}})\\ \nonumber
&=(\phi-s\phi')\bar{\textbf{n}}_{a}+\phi''\beta(\bar{\textbf{n}}_{a})(W-s\bar{\textbf{n}}),
\end{align}
where $\textbf{n}_{a}$ is the derivative of $\textbf{n}$ with respect to $u$ or $v$.
Let ~$W^{T}=W-s\bar{\textbf{n}}=W^{a}\textbf r_{a}$, by calculation, we can obtain
\begin{equation}\label{10}
W^{u}=<W^{T}, \textbf r_{u}>=\beta(\textbf r_{u})=b^{1}\cos v+b^{2}\sin v,
\end{equation}
\begin{equation}\label{11}
W^{v}=\frac{\beta(\textbf r_{v})}{G}=\frac{1}{G}(-ub^{1}\sin v+ub^{2}\cos v+ab^{3}).
\end{equation}
Then by \eqref{8},
\begin{align*}
\mathcal{A}(\textbf{r}_{u})&=-\textbf{n}_{u}=-(\phi-s\phi')\bar{\textbf{n}}_{u}-\phi''\beta(\bar{\textbf{n}}_{u})(-s\bar{\textbf{n}}+W)\\
&=-(\phi-s\phi')\mu_{1}\textbf{r}_{v}-\phi''\beta(\mu_{1}\textbf{r}_{v})W^{T}.
\end{align*}

Let $\phi-\beta\phi'=\varphi$,
then the above formula is equal to
\begin{align*}
\mathcal{A}(\textbf r_{u})
&=-\varphi \mu_{1}\textbf r_{v}-\phi''\beta(\mu_{1}\textbf r_{v})W^{T}\\
&=-\varphi \mu_{1}\textbf r_{v}-\mu_{1}\phi''\beta(\textbf r_{v})(\beta(\textbf r_{u})\textbf r_{u}+\frac{\beta(\textbf r_{v})}{G}\textbf r_{v}).
\end{align*}
Similarly, we have
\begin{align*}
\mathcal {A}(\textbf r_{v})&=-\textbf{n}_{v}=-(\phi-\beta\phi')\mu_{2}\textbf r_{u}-\phi''\beta(\bar{\textbf{n}}_{v})W^{T}\\
&=-\varphi \mu_{2}\textbf r_{u}-\mu_{2}\phi''\beta(\textbf r_{u})(\beta(\textbf r_{u})\textbf r_{u}+\frac{\beta(\textbf r_{v})}{G}\textbf r_{v}).
\end{align*}
So
$$\mathcal{A}\left(
                      \begin{array}{cc}
                        \textbf r_{u}  \\
                        \textbf r_{v}  \\
                      \end{array}
                    \right)=-\left(
                      \begin{array}{cc}
                        \mu_{1}\phi''\beta(\textbf r_{v})\beta(\textbf r_{u}) &\mu_{1}(\varphi+\frac{\phi''\beta^{2}(\textbf r_{v})}{G})\\
                        \mu_{2}(\varphi+\phi''\beta^{2}(\textbf r_{u})) & \frac{\mu_{2}}{G}\phi''\beta(\textbf r_{v})\beta(\textbf r_{u}) \\
                      \end{array}
                    \right)\left(
                      \begin{array}{cc}
                        \textbf r_{u}  \\
                        \textbf r_{v}  \\
                      \end{array}
                    \right)=\omega\left(
                      \begin{array}{cc}
                        \textbf r_{u}  \\
                        \textbf r_{v}  \\
                      \end{array}
                    \right).
$$
If set $W=(0,0,b)$, then $W^{u}=0, W^{v}=\frac{ab}{G}.$ By \eqref{10}, \eqref{11}, we get
$\beta(\textbf r_{u})=0,\beta(\textbf r_{v})=ab.$ So
the coefficient matrix becomes
$$\omega=\left(
                      \begin{array}{cc}
                        0 & -\mu_{1}(\varphi+\frac{1}{G}\phi''a^{2}b^{2}) \\
                        - \mu_{2}\varphi & 0 \\
                      \end{array}
                    \right).
$$
 The characteristic polynomial of the matrix is
$$  \lambda^{2}=\mu_{1}\mu_{2}\varphi (\varphi+\frac{1}{G}\phi''a^{2}b^{2}).$$
If
\begin{equation}\label{12}
\lambda^{2}=1,
\end{equation}
Because $\mu_{2}=G\mu_{1}$, we can get
\begin{align*}
1=&G\mu_{1}^{2}\varphi(\varphi+\frac{1}{G}\phi''a^{2}b^{2})=G\mu_{1}^{2}\varphi^{2}+\mu_{1}^{2}\varphi\phi''a^{2}b^{2}\\
&=\varphi\mu_{1}^{2}(G\varphi+\phi''a^{2}b^{2})=\frac{a^{2}\varphi}{G^{3}}(G\varphi+\phi''a^{2}b^{2}).
\end{align*}
Since $s= \beta(\bar{\textbf{n}})=\frac{bu}{\sqrt{u^2+a^2}}$,~$G=u^2+a^2=\frac{a^2b^2}{b^2-\beta^2},$~the above formula holds if and only if
$$2\varphi'\varphi-\frac{2s}{b^2-s^2}\varphi^2=\frac{-2sa^2b^4}{(b^2-s^2)^3}.$$
Let $f=\varphi ^2,$ we can obtain
\begin{equation}\label{13}
f'-\frac{2s}{b^2-s^2}f+\frac{2a^2b^4s}{(b^2-s^2)^3}=0.
\end{equation}
Set $f(0)=1$, by solving this equation, we can get \\
$$f=\frac{b^2(b^2-a^2s^2-s^2)}{(b^2-s^2)^2},$$
so
$$\varphi =\frac{b\sqrt {b^2-a^2s^2-s^2}}{b^2-s^2}.$$
By solving
\begin{equation}\label{100}
\phi-\beta\phi'=\varphi,
\end{equation}
we can obtain
\begin{equation}\label{14}
\phi=|s|\left(\int\frac{-b\sqrt{(b^2-a^2s^2-s^2)}}{(b^2-s^2)s|s|}ds\right).
\end{equation}
If $\phi$ is defined by
\begin{equation}\label{101}
\phi=
\begin{cases}
s \int_{s}^{c}\frac{b\sqrt{(b^2-a^2t^2-t^2)}}{(b^2-t^2)t^2}dt, & {s>0},\\
-s \int_{-c}^{s}\frac{b\sqrt{(b^2-a^2t^2-t^2)}}{(b^2-t^2)t^2}dt,  &{s<0},\\
   \end{cases}
\end{equation}
 where $c=\frac{b}{\sqrt{1+a^2}}$.
By direct calculation, we can get
\begin{equation}\label{183}\int\frac{-b\sqrt{b^2-(1+a^2)s^2}}{(b^2-s^2)s^2}ds
=\frac{\sqrt{b^2-(1+a^2)s^2}}{bs}-\frac{a}{b}\arctan\frac{\sqrt{b^2-(1+a^2)s^2}}{as}.
\end{equation}
So
\begin{equation}\label{139}
\phi=|s|\Big(\frac{\sqrt{b^2-(1+a^2)s^2}}{b|s|}-\frac{a}{b}\arctan\frac{\sqrt{b^2-(1+a^2)s^2}}{a|s|}\Big),~{0<|s|<c},
\end{equation}
 then in $(-c,0)\cup(0,c),~\phi(s)$ is a smooth positive function and
\begin{equation}
 \phi-s\phi'+(b^2-s^2)\phi''=\varphi-(b^2-s^2)\frac{\varphi'}{s}.
 \end{equation}
Because $\varphi >0$,
\begin{align*}
\varphi( \varphi-(b^2-s^2)\frac{\varphi'}{s})
&=f-\frac{(b^2-s^2)}{2s}f'\\
&=f-\frac{(b^2-s^2)}{2s}(\frac{2s}{(b^2-s^2)}f-\frac{2a^2b^4s}{(b^2-s^2)^3)})\\&=\frac{a^2b^4}{(b^2-s^2)^2}>0,
 \end{align*}
thus $$\phi-s\phi'+(b^2-s^2)\phi''>0.$$
$F^\ast$ satisfies the conditions of the Example \ref{3.0} in $(-c,0)\cup(0,c)$, so $F^\ast$ is a conic Minkowski-$(\alpha,\beta)$ metric, which can be expressed as$$F^{\ast}=\sqrt{\xi_{1}^2+\xi_{2}^2-a^{2}\xi_{3}^{2}}-a\xi_{3}\arctan\frac{\sqrt{\xi_{1}^2+\xi_{2}^2-a^{2}\xi_{3}^{2}}}{a\xi_{3}},$$ whose conic domain is $$A^\ast =\{ \xi\in \mathbb{R}^3 |\xi_1^2+\xi_2^2>a^2\xi_3^2 ,\xi_3^2>0\}.$$ Thus $F$ is a conic Minkowski-$(\alpha,\beta)$ metric. We can prove that $A^\ast$ is a applicable cone. Suppose~$\mathcal{L}_{F^\ast}(A^\ast)=A\subset \mathbb{R}^3$,~we only need to prove~$\mathcal{L}_{F^\ast} :A^\ast\rightarrow A$~is a bijection.

For simplicity, set~$(x,y,z)=\xi=(\xi_{1},\xi_{2},\xi_{3})$. Surface~$\Sigma_{\pm}:\{~\xi\in A^\ast|F^\ast(\xi)=1\}$~is obtained by rotating~$L_{\pm} :\{~\xi\in A^\ast|F^\ast(\xi)=1,~y=0,~x>0\}$ around the $z$ axis. Where $L_{+}$~represents the curve when $s > 0$, and $L_{-}$~represents the curve when $s < 0$. Therefore, it is only necessary to consider whether~$L_{\pm}$ and $\Sigma_{\pm}$ have common tangent.

Let~$\xi$~be a point of~$L_{\pm}$,~then $F^\ast=\alpha^\ast\phi(s),~\alpha^\ast=|\xi|=\sqrt{x^2+z^2}$,~that is, $s=\frac{\beta^\ast}{\alpha^\ast}=\frac{bz}{|\xi|}.$~
Set~$\rho=|\xi|$, the angle between $\xi$ and $x$ axis is~$\theta$.~Because~$F^\ast(\xi)=1,$~we can obtain $\rho=\frac{1}{\phi(s)},$~then~$L_{\pm}: \textbf{r}(s)=(x(s),z(s))$~can be expressed as \begin{equation}\begin{cases} x(s)=\rho(s)\cos\theta,\\
z(s)=\rho(s)\sin\theta.
\end{cases}
\end{equation}
Since~$\frac{\xi_{3}}{|\xi|}=\sin\theta$,~$s=b\sin\theta,$~then we have
\begin{equation}\label{123}\begin{cases} x(s)=\frac{\sqrt{b^2-s^2}}{b\phi},\\
z(s)=\frac{s}{b\phi}.
\end{cases}
\end{equation}
By\eqref{139}~and~\eqref{123},~we can get~$\lim\limits_{s\rightarrow 0}\textbf{r}(s)=(1,0)$~and
~$$ x'(s)=\frac{-s(\phi-s\phi')-b^2\phi'}{\phi^{2}b\sqrt{b^2-s^2}},~z'(s)=\frac{\phi-s\phi'}{b\phi^2}.$$
Because~$z'(s)>0$,~$L_{\pm}$~can be expressed as~$x=x(z).$~By direct calculation, we can get
\begin{align}
\frac{dx}{dz}&=\frac{-s(\phi-s\phi')-b^2\phi'}{(\phi-s\phi')\sqrt{b^2-s^2}}\nonumber\\
&=\frac{1}{\sqrt{b^2-s^2}}(-s-b^2\frac{\phi'}{\phi-s\phi'})\nonumber\\
&=a\frac{\sqrt{b^2-s^2}}{\sqrt{b^2-(1+a^2)s^2}}\arctan\frac{\sqrt{b^2-(1+a^2)s^2}}{as}.\label{113}
\end{align}
According to \eqref{113},~we can obtain
$$\frac{dx}{dz}>0,~\forall s\in(0,c),~~~~\frac{dx}{dz}<0,~\forall s\in(-c,0).~$$
Furthermore, the tangent of~$L_{\pm}$~at any point is
$$\frac{x-\frac{\sqrt{b^2-s^2}}{b\phi}}{-\Big(s(\phi-s\phi')+b^2\phi'\Big)}=\frac{z-\frac{s}{b\phi}}{(\phi-s\phi')\sqrt{b^2-s^2}}.$$
When $z=0$, by calculation, we can get $x=\frac{b}{(\phi-s\phi')\sqrt{b^2-s^2}}>0$. So~$L_{+}$~and surface~$\Sigma_{\pm}$~have no common tangent.~Similarly,~$L_{-}$~and surface~$\Sigma_{\pm}$~have no common tangent. Therefore~$\mathcal{L}_{F^\ast}$~is a injection, thus, the dual metric~$F$~of~$F^{\ast}$~is well defined, and is also a conic~Minkowski-$(\alpha,\beta)$~metric. 

Due to $$\lim_{|s|\rightarrow c^{-}}\frac{dx}{dz}=\pm a,~~~~\lim_{s\rightarrow 0^{\pm }}\frac{dx}{dz}=\pm \frac{\pi a}{2},$$and~$\mathcal{L}_{F^\ast}(\xi)$~is the normal vector of ~$L_{\pm}$,~$\mathcal{L}_{F^\ast}(\xi)$~is parallel to~$(1,-\frac{dx}{dz})$.~So the conic domain of $F$ is
$$ A_F=\Big\{y\in \mathbb{R}^3~|~0<\frac{4}{\pi^2a^2}(y^{3})^2<(y^{1})^2+(y^{2})^2<\frac{1}{a^2}(y^{3})^{2}\Big\}.$$

Since the above process is reversible, we can obtain \eqref{12}. In order to make $\textbf{n}\in AM$, we need to let $\nu \in A^\ast M$, which means $\bar{\textbf{n}}=\bar \nu \in A^\ast M$.~Therefore, when $0<u<1$, $\textbf r=(u\cos v,u\sin v,av)$ is a local isoparametric hypersurface with constant principal curvatures $\pm 1$.
\endproof
\section{Isoparametric Hypersurfaces in  Kropina Spaces}
\subsection{Conic submanifolds in  Kropina Spaces}
A Finsler metric $F$ is a Randers metric if and only if it is the solution of Zermelo navigation problem on a Riemannian space $(M, h)$ under the influence of a force field $W$ with $\|W\|_{h}<1$, where $\|W\|_{h}$ denotes the length of $W$ with respect to Riemannian metric $h$.

Similarly, Kropina metrics can also be characterized as the solution of the Zermelo navigation problem on a Riemannian space $(M, h)$ under the influence of a force field $W$ with $\|W\|_{h}=1$. Concretely, assuming that $h=\sqrt{ h_{ij}(x)y^{i}y^{j}}$  and $W =W^{i}\frac{\partial}{\partial x^i} $,
then the solution of the Zermelo navigation problem is a Kropina metric given by
\begin{equation}
F=\frac{h^2}{2W_0}.
\end{equation}
For a Kropina metric $F$, at each $x\in M$, the conic domain of $T_x M$ is defined as following
\begin{center}
$A_xM:=\{y=y^i\frac{\partial}{\partial y^i}\in T_x M\vert W_0(x,y)=W_i(x)y^i>0\}$,\end{center}
where its boundary is the hyperplane $\{y=y^i\frac{\partial}{\partial y^i }\in T_x M\vert W_i(x)y^i=0\}$.~Since the set $\{y\in A_xM|F(y)<1\}$ is a parallel shifting of $\{y\in T_xM|h(y)<1\}$, the \emph {BH volume form} of Kropina metrics can be well defined.

Denote the dual metric of $h$ by $h^\ast$. Then the dual metric of $F$ can be expressed as\cite{DH1}
\begin{equation}\label{2.3}
F^\ast=h^\ast+W^0=\sqrt{h^{ij}\xi_i\xi_j}+W^i\xi_i,~~~~~~~\xi=\xi_idx^i\in A_x^\ast M.
\end{equation}
Take $\textbf{n}\in \mathcal{N}^{0}(N)$. From \eqref{2.3}, we know that
$$\textbf{n}^i=F^{\ast}_{\xi_i}(\nu)=\frac{h^{ij}\nu_j}{h^\ast(\nu)}+W^i.$$
Denote $\bar{\textbf{n}}=\frac{h^{ij}\nu_j}{h^\ast(\nu)}$. Then $\bf{\bar n}$ is a unit normal vector filed of $N$ with respect to $h$. Thus
\begin{equation}\label{1.4}
\textbf{n}=\bar{\textbf{n}}+W.
\end{equation}
$N$ is a conic submanifold of $(M, F)$ if and only if $\textbf{n}\in AM$, that is, $W_{0}({\textbf{n}})=W_{0}(\bar{\textbf{n}})+1=\langle\bar{\textbf{n}}, W\rangle_h+1>0$.
\begin{lem}
Let $\Phi: N \rightarrow (M, F) $ be a conic submanifold in a Kropina space $(M, F)$ with the navigation data $(h, W)$, then for any smooth section $\textbf{n}$ of $\mathcal{N}^{0}(N)$, $\bf{\bar n}$$=\textbf{n}- W$ is a unit normal vector filed of $N$ with respect to $h$  satisfying $\bar{\textbf{n}}\neq-W$, and the induced metric $\hat g_{\textbf{n} }=\Phi^\ast g_{\textbf{n}}$ is conformal to $\bar h = \Phi^\ast h$ and satisfies
$$\hat g_{\textbf{n} }=\frac{1}{W_{0}(\bf{n})}\bar h=\frac{1}{W_{0}(\bar{\textbf{n}})+1}\bar h.$$
\end{lem}
\proof
By direct computation, we can obtain
\begin{equation}\label{1.4.1} g_{ij}=\frac{F}{W_{0}}(h_{ij}-\frac{2}{W_0}hh_{y^{i}}W_j+\frac{h^2W_{i}W_{j}}{W_{0}^2})+F_{y^i}F_{y^j}.\end{equation}
Let $(u^a) = (u^1, . . . , u^n)$ be the local coordinates on $N$ and $d\Phi  = \Phi_a^i du^{a}\otimes \frac{\partial}{\partial x^i}$ . Then\\
\begin{align*}
F_{y^i}(\bf{n}) &=\frac{2h(\textbf{n})h_{y^i}(\textbf{n})W_{0}(\textbf{n})-W_{i}h^{2}(\bf{n})}{2W_{0}^2(\bf{n})}\\
&=\frac{h(\textbf{n})}{W_{0}(\textbf{n})}h_{y^i}(\textbf{n})-W_{i}\frac{h^2(\textbf{n})}{2W_{0}^2(\textbf{n})}.
\end{align*}
Following from $F(\textbf{n})=1$, we can obtain
$2W_{0}(\textbf{n})=h^2(\textbf{n})$ and
\begin{align}\label{1.4.2}
h_{y^i}(\textbf{n})\Phi_a^i&= \left(F_{y^i}(\textbf{n})\Phi_a^i+\frac{h^2(\textbf{n})}{2W_{0}^2(\textbf{n})}W_{i}\Phi_a^i\right)\frac{W_{0}(\textbf{n})}{h(\textbf{n})}=\frac{1}{h(\textbf{n})}W_{i}\Phi_a^i.
\end{align}
So
\begin{align*}
(\hat g_\textbf{n})_{ab}&=g_{ij}(\textbf{n})\Phi_a^i\Phi_b^j\\
&=\left[\frac{F(\textbf{n})}{W_{0}(\textbf{n})}\left(h_{ij}-\frac{2}{W_{0}(\textbf{n})}h(\textbf{n})h_{y^i}(\textbf{n})W_{j}+\frac{h^2(\textbf{n})W_{i}W_{j}}{W_{0}^2(\textbf{n})}\right)+F_{y^i}(\textbf{n})F_{y^j}(\textbf{n})\right]\Phi_a^i\Phi_b^j\\
&=\frac{1}{W_{0}(\textbf{n})}\left[h_{ij}-\frac{2}{W_{0}(\textbf{n})}h(\textbf{n})h_{y^i}(\textbf{n})W_{j}+\frac{2W_{i}W_{j}}{W_{0}(\textbf{n})}\right]\Phi_a^i\Phi_b^j\\
&=\frac{1}{W_{0}(\textbf{n})}\left(\bar h_{ab}-\frac{2}{W_{0}(\textbf{n})}h(\textbf{n})\frac{1}{h(\textbf{n})}W_{i}\Phi_a^i\Phi_b^jW_{j}+\frac{2W_{i}W_{j}}{W_{0}(\textbf{n})}\Phi_a^i\Phi_b^j\right)\\
&=\frac{1}{W_{0}(\textbf{n})}\bar h_{ab},
\end{align*}
where
$\bar h_{ab}=h_{ij}\Phi_a^i\Phi_b^j.$
\endproof
Denote
$$r_{ij}=\frac{1}{2}(W_{i|j}+ W_{j|i}),~~~s_{ij}=\frac{1}{2}(W_{i|j}- W_{j|i}),$$
$$r_{j}=W^{i}r_{ij},~~~r=r_{j}W^{j},~~~r^{i}=h^{ik}r_{k},~~~r^{i}_{~j}=h^{ik}r_{kj},$$
$$r_{i0}=r_{ij}y^{j},~~~r^{i}_{~0}=r^{i}_{~j}y^{j},~~r_{0}=r_{j}y^{j},~~r_{00}=r_{ij}y^{i}y^{j},$$
$$s_{j}=W^{i}s_{ij},~~~s_{i0}=s_{ij}y^{j},~~~~s_{0}=s_{i}y^{i},~~~s_{0}^{i}=s_{j}^{i}y^{j},$$
where $|$ denotes the covariant differentiation with respect to $h$.
If $W$ is a Killing vector field of constant length $\|W\|_{h}=1$ on $M$, then $r_{ij}=0,~s_{i}=0.$ From\cite{ZS1}, we get
$$G^i=\bar G^i -Fs_{0}^{i},$$
where $G^i$ and $\bar G^i$ are the geodesic coefficients of $F$ and $h$, respectively.

\begin{lem}
Let $\Phi: N\rightarrow (M,F)$ be a conic submanifold in a Kropina space $(M,F)$ with the navigation data $(h,W)$, where $W$ is a Killing vector field of constant length $\|W\|_{h}=1$ on $M$, then for any $\textbf{n}\in \mathcal{N}^{0}(N)$ and $X\in TN$,
\begin{equation}\label{31}
D_{X}^{ \textbf{n}}\textbf n=\nabla _{X}^{h}{\bar{\textbf{n}}}.
\end{equation}
\end{lem}
\proof
\begin{align*}
D_{X}^{\textbf n}\textbf n&=(\textbf n^{i}_{x^{j}}+N^i_{j}(\textbf n))\Phi_{a}^{j} X^a\frac{\partial}{\partial x^i}\\
                     & =(\textbf n^{i}_{x^{j}}+ \bar N^i_{j}(\textbf n)-F_{y^j}(\textbf n)s_{0}^{i}-F(\textbf n)s_j^i) \Phi_{a}^{j} X^a\frac{\partial}{\partial x^i}\\
                        &=(\textbf n^{i}_{x^{j}}+ \bar N^i_{j}(\textbf n)-W^i_{|j}  )\Phi_{a}^{j} X^a\frac{\partial}{\partial x^i} \\
                        &=\nabla _{X}^{h} \textbf n-W^i_{|j}\Phi_{a}^{j} X^a\frac{\partial}{\partial x^i} \\
                        &=\nabla _{X}^{h} ({\bf\bar n}+W)-W^i_{|j}\Phi_{a}^{j} X^a\frac{\partial}{\partial x^i} \\
                        & =\nabla _{X}^{h} {\bf\bar n}+\nabla _{X}^{h} W-W^i_{|j}\Phi_{a}^{j} X^a\frac{\partial}{\partial x^i} \\
                       & =\nabla _{X}^{h} \bf{\bar n}.
\end{align*}
\endproof
\begin{remark}
The condition that $W$ is a Killing vector field can be changed to that $F$ has isotropic $\bf S$-curvature $S=(n+1)k(x)F$, which is mentioned in\cite{HD}.~But in this case, according to\cite{X}, $k(x)\equiv0$, so $W$ becomes a Killing vector field automatically.
\end{remark}
\begin{thm}\label{001}
Let $N$ be a conic submanifold in a Kropina space $(M, F, d\mu_{BH})$ with the navigation data $(h, W)$, where $W$ is a Killing vector field of constant length $\|W\|_{h}=1$ on $M$. For any $\textbf{n}=\bar{\textbf{n}}+W \in \mathcal{N}^{0}(N),$ the shape operators of $N$ in Kropina space $(M, F)$ and Riemannian space $(M, h)$, $\mathcal{A}_{\textbf n}$ and $\mathcal {\bar A}_{\bar{\textbf{n}}}$, have the same principal vectors and principal curvatures.
\end{thm}
\proof
Set $X=X^a\frac{\partial}{\partial u^a}$ and $\Phi_a=d\Phi \left(\frac{\partial}{\partial u^a}\right)$.
By \eqref{30}, \eqref{1.4.1}, \eqref{1.4.2} and \eqref{31}, we know that
\begin{align*}
\mathcal{A}_{\textbf n}X=-[D_{X}^{\textbf n}\textbf n]_{g_\textbf n}^T
&=-g_{\textbf{n}}(\nabla _{X}^{h} {\bf \bar n},\Phi_a)(\hat g_{\textbf{n}})^{ab}\frac{\partial}{\partial u^b}\\
&=-\frac{1}{W_{0}\bf (n)}\left( h_{ij}-\frac{2}{W_{0}(\textbf{n})}h(\textbf{n})h_{y^i}(\textbf{n})W_{j}+\frac{2W_{i}W_{j}}{W_{0}(\textbf{n})}\right)
\Phi_a^i X^c {\bf\bar n}^j_{|c}(\hat{g}_\textbf n)^{ab}\frac{\partial}{\partial u^b}\\
&=-\frac{1}{W_{0}\bf (n)}h_{ij}\Phi_a^i X^c{\bf\bar n}^j_{|c}{W_{0}\bf (n)}\bar h^{ab}\frac{\partial}{\partial u^b}\\
&+\left(\frac{2h(\textbf{n})}{W_{0}^2({\textbf{n}})}h_{y^i}(\textbf{n})\Phi_a^i W_{j}-\frac{2W_{i}W_{j}}{W_{0}^2(\textbf{n})}\Phi_a^i\right)X^c {\bf\bar n}^j_{|c}(\hat{g}_\textbf n)^{ab}\frac{\partial}{\partial u^b}\\
&=-h_{ij}\Phi_a^i X^c{\bf\bar n}^j_{|c}\bar h^{ab}\frac{\partial}{\partial u^b}\\
&=-[\nabla _{X}^{h} {\bf\bar n} ]_{g_\textbf n}^T\\
&=\mathcal{\bar A}_{\bf\bar n}X.
\end{align*}
Thus $\mathcal{A}_{\textbf n}$ and $\mathcal{\bar A}_{\bf\bar n}$ have the same principal vectors and  principal curvatures.
\endproof
The following corollaries are immediate consequences of Theorem \ref{001}, so we skip their proofs.
\begin{corr}
In a Kropina space $(M, F, d\mu_{BH})$ with the navigation data $(h, W)$, where $W$ is a Killing vector field of constant length $\|W\|_{h}=1$, a conic submanifold $N$ is totally umbilic if and only if it is totally umbilic in Riemannian space $(M, h)$.
\end{corr}
\begin{corr}
In a Kropina space  $(M, F, d\mu_{BH})$ with the navigation data $(h, W)$, where $W$ is a Killing vector field of constant length $\|W\|_{h}=1$, the principal curvatures of a conic submanifold are all constant if and only if its principal curvatures in Riemannian space $(M, h)$ are all constant.
\end{corr}
\begin{corr}
In a Kropina space $(M, F, d\mu_{BH})$ with the navigation data $(h, W)$, where $W$ is a Killing vector field of constant length $\|W\|_{h}=1$, a conic submanifold $N$ has constant mean curvature if and only if $N$ also has constant mean curvature in Riemannian space $(M, h)$. Especially, $N$ is minimal if and only if it is minimal in Riemannian space $(M, h)$.
\end{corr}
\subsection{Proof of Theorem 1.3}

\begin{lem}\cite{YRSS}\label{131}
Let $(M, F)$ be an $m(\geqslant2)$ dimensional Kropina space, $h=\sqrt{ h_{ij}(x)y^{i}y^{j}}$ and a vector field
$W =W^{i}\frac{\partial}{\partial x^i} $ of constant length $|| W||_{h}=1$ on $M$,
then the Kropina space $(M, F)$ is of constant curvature $K$ if and only if the following conditions hold:\\
(1) $W_{i|j} + W_{j|i} = 0$, that is, $W =W^{i}\frac{\partial}{\partial x^i}$ is a Killing vector field. \\
(2) The Riemannian space $(M, h)$ is of constant sectional curvature $K$.
\end{lem}
\begin {lem}\cite{YRSS}\label{135}
The only manifolds (up to Riemannian local isometry) that admits \emph{CC Kropina structures}($W$ is a unit length Killing vector field on the Riemannian space $(M, h)$ of constant sectional curvature)are the Euclidean space $E^{m}$, $m\geqslant2$ and odd dimensional spheres $S^{2m-1}$, $m\geqslant 2$.
\end{lem}
\proof
The first half of Theorem \ref{130} follows from Lemma \ref{131}, Theorem \ref{132} and Theorem \ref{001}. Moreover, from Lemma \ref{135} and references             \cite{C38} -\cite{M13}, we can give the complete classifications of isoparametric hypersurfaces in a Kropina space with constant flag curvature, and the results are summarized in the following table.
\endproof

\begin{table}[h]
  \begin{center}
  \caption{Classifications for isoparametric hypersurfaces in $(M(K), F)$}
    \begin{tabular}{|c|c|c|c|c|c|l|c|}
  \hline
     $K_F=c$ &$\mathbf{S}$-curv. & $M(c)$ &$g$ &dim$N$ & mul.& $N$ is an open subset of&main \\
      &&&&&& following hypersurfaces &ref. \\
  \cline{1-8}
      \multirow{2}{*}{$K=0$} &\multirow{2}{*}{$S=0$}  & $\mathbb{R}^{m}$ &$g$=1& $m - 1$&$m - 1$ & a hypersphere $\mathbb{S}^{m - 1}$&\multirow{4}{*}{}\\
      &&$\|W\|_h=1$&&&& or a hyperplane $\mathbb{R}^{m - 1}$&\cite{C38}\\
 \cline{4-7}
       & &&$g$=2& $m - 1$&($n$,$m - n - 1$) &a cylinder $\mathbb{S}^{n}\times$ $\mathbb{R}^{m - n - 1}$ &\cite{SF}\\
  \cline{1-8}
       \multirow{6}{*}{\emph{$K=1$}}& \multirow{6}{*}{\emph{$S=0$}}&\multirow{6}{*}{$\mathbb{S}^{2m - 1}$}& $g$=1 &$2m- 2$  &$2m- 2$ & a great or small &\\
       &&&&&& hypersphere &\\
  \cline{4-7}
      &&& $g$=2&$2m - 2$&($n$,$2m - n - 2$)&a Clifford torus &\cite{C39}\\
      &&&&&&$S^{n}(r)\times S^{2m - n - 2}(s)$, &\cite{C40}\\
      &&&&&& $r^{2}+s^{2}=1$&\\
  \cline{4-7}
  &&&$g=3$&6&(2,2)&a tube over a standard&\\
    \cline{5-6}
       &&&&12&(4,4)& Veronese embedding of&\\
  \cline{5-6}
   &&&&24&(8,8)& $\mathbb{F}$P into $S^{3n+1}$, where&\\
   &&&&&& $\mathbb{F}$=$\mathbb{C}$,$\mathbb{H}$ or $\mathbb{O}$, for &\\
      &&&&&& $n=2,4,8$, respectively&\\
 \cline{4-8}
      &&$\|W\|_h=1$&\multirow{4}{*}{$g$=4}&$2(n_{1}+n_{2})$&$(n_{1},n_{2})$&\multirow{3}{*}{OT-FKM type or}  &\cite{CS1}\\
      &&&&$n_{2}\geq2n_{1}-1$&&&\cite{FKM81}\\
  \cline{5-6}
      &&&&8&(2,2)&homogeneous &\cite{M80}\\
  \cline{5-6}
      &&&&18&(4,5)&&\cite{OT76}\\
  \cline{5-6}
      &&&&14&(3,4)&&\cite{C2}\\
      \cline{5-6}
      &&&&30&(6,9)&&\cite{C3}\\
  \cline{5-6}
      &&&&30&(7,8)&&\cite{C1}\\
  \cline{4-8}
   &&&\multirow{2}{*}{$g$=6}&6&(1,1)&\multirow{2}{*}{homogeneous} &\cite{DN85}\\
  \cline{5-6}
      &&&&12&(2,2)& &\cite{M13}\\
      \hline
    \end{tabular}
  \end{center}
\end{table}
\endproof

\noindent {Qun He;~~~Xin Huang} \\
\noindent {School of Mathematics Science, Tongji University,
 Shanghai, 200092, China.} \\
E-mail address:  hequn$@$tongji.edu.cn; ~~~~~1930913$@$tongji.edu.cn

\noindent {Peilong Dong}\\
\noindent {School of Mathematics and Statistics, Zhengzhou Normal University,\\
Zhengzhou Henan 450044, China.}\\
E-mail: dpl2021@163.com
\end{document}